\pdfoutput=1
\RequirePackage{ifpdf}
\ifpdf %We are running pdfTeX in pdf mode
\documentclass[pdftex]{sigma}
\else
\documentclass{sigma}
\fi

\usepackage[all]{xy}
\numberwithin{theorem}{section}
\numberwithin{proposition}{section}
\numberwithin{lemma}{section}
\numberwithin{corollary}{section}
\numberwithin{definition}{section}
\numberwithin{example}{section}
\numberwithin{remark}{section}
\numberwithin{note}{section}

\DeclareMathOperator{\Ker}{Ker}
\DeclareMathOperator{\End}{End}

\numberwithin{equation}{section}

\begin{document}

\allowdisplaybreaks

\renewcommand{\thefootnote}{$\star$}

\renewcommand{\PaperNumber}{040}

\FirstPageHeading

\ShortArticleName{From Twisted Quantum Loop Algebras to Twisted Yangians}

\ArticleName{From Twisted Quantum Loop Algebras\\
to Twisted Yangians\footnote{This paper is a~contribution to the Special Issue on Exact Solvability and Symmetry Avatars
in honour of Luc Vinet.
The full collection is available at
\href{http://www.emis.de/journals/SIGMA/ESSA2014.html}{http://www.emis.de/journals/SIGMA/ESSA2014.html}}}

\Author{Patrick CONNER~$^\dag$ and Nicolas GUAY~$^\ddag$}

\AuthorNameForHeading{P.~Conner and N.~Guay}

\Address{$^\dag$~Science Department, Red Deer College, Red Deer, Alberta T4N 5H5, Canada}
\EmailD{\href{mailto:patrick.conner@rdc.ab.ca}{patrick.conner@rdc.ab.ca}}

\Address{$^\ddag$~Department of Mathematical and Statistical Sciences, University of Alberta,\\
\hphantom{$^\ddag$}~CAB 632, Edmonton, Alberta T6G 2G1, Canada}
\EmailD{\href{mailto:nguay@ualberta.ca}{nguay@ualberta.ca}}
\URLaddressD{\url{http://www.ualberta.ca/~nguay/}}

\ArticleDates{Received January 30, 2015, in f\/inal form May 01, 2015; Published online May 14, 2015}

\Abstract{We prove how the Yangian of $\mathfrak{gl}_N$ in its RTT presentation and Olshanski's twisted
Yangians for the orthogonal and symplectic Lie algebras can be obtained by a~degeneration process from the corresponding
quantum loop algebra and some of its twisted analogues.}

\Keywords{twisted Yangians; twisted quantum loop algebras; degeneration; RTT-pre\-sen\-ta\-tion}

\Classification{17B37; 16T20; 81R50}

\renewcommand{\thefootnote}{\arabic{footnote}}
\setcounter{footnote}{0}

\section{Introduction}

Yangians are quantum groups which appeared over thirty years ago.
Since that time, it has been shown that they play a~role in the theory of integrable systems by controlling some of
their symmetries~\cite{Be,BMM,ChPr,HHTBP,Ma2,MRSZ}.
It has been discovered that, for integrable systems with boundary, symmetries are better captured by twisted Yangians,
which are coideal subalgebras of Yangians~\cite{dRGMa,DMS, MaRe1,MaRe2, Ma1}.
This has provided motivation to mathematicians to study the representation theory of Yangians and twisted Yangians:
see~\cite{Mo} and references therein, in particular the work of A.~Molev and M.~Nazarov.

It has been well known for a~long time that quantum loop algebras and Yangians associated to complex semisimple Lie
algebras have similar representation theories although, until very recently, explicit connections between their module
categories (other than just analogies and similarities) were not available.
The papers~\cite{GaTL2,GaTL3} seem to be the f\/irst ones to establish rigorously important equivalences between certain
categories of modules for both of those quantum groups.

It is also commonly known that Yangians are limit forms of quantum loop algebras: this was explicitly stated
in~\cite{Dr1}, but a~proof appeared only much later in~\cite{GuMa} although such a~proof was certainly known to some
experts.
(One early result in that direction can be found in Appendix C in~\cite{BeLe}.) Stronger results were obtained at about
the same time by V.~Toledano Laredo and S.~Gautam: see~\cite{GaTL1}.
The results of these two papers apply to Yangians associated to complex semisimple Lie algebras presented using what is
commonly called Drinfeld's second realization~\cite{Dr2}.
Our f\/irst goal in this article is to prove that the Yangian of $\mathfrak{gl}_N$ is isomorphic to a~limit
form $\widetilde{Y}(\mathfrak{gl}_N)$ of the quantum loop algebra of $\mathfrak{gl}_N$ using the
generators from the RTT-presentation.
(Some indications about how to show this are given in~\cite{Ch}, but we follow a~dif\/ferent approach.) It is not clear if
this could be used to deduce the main theorem of~\cite{GuMa} in the \mbox{$\mathfrak{s}\mathfrak{l}_n$-case} or vice-versa:
this would probably require checking that the isomorphisms in the main theorem of~\cite{GuMa} and in
Theorem~\ref{mainthm} below are compatible with the isomorphisms between the two realizations of quantum loop algebras
and Yangians for $\mathfrak{gl}_n$ established in~\cite{DiFr} and~\cite{BrKl}.

In the second half of~\cite{GuMa}, it was established that certain degenerate forms of twisted quantum loop algebras
associated to a~complex semisimple Lie algebra $\mathfrak{g}$ and a~Dynkin diagram automorphism of $\mathfrak{g}$ are
isomorphic to the ordinary Yangian of $\mathfrak{g}$.
This explains partially why there is no known twisted Yangian attached to such data in general.
The two cases considered in section~\ref{twsec} below are however dif\/ferent: the twisted quantum loop algebras of type
AI and AII are associated to an involution~$\theta$ of $\mathfrak{gl}_N$ and are quantizations of the
enveloping algebra of the twisted loop algebra which is the f\/ixed point of the involution~$\Theta$ on
$\mathfrak{gl}_N\otimes_{\mathbb{C}} \mathbb{C}[s,s^{-1}]$ given by $\Theta(X \otimes s) = \theta(X) \otimes
s^{-1}$, not $\Theta(X \otimes s) =-\theta(X) \otimes s$ as would be the case with a~twisted af\/f\/ine Kac--Moody algebra
associated to a~Dynkin diagram automorphism~$\theta$ of order~$2$.
Theorem~\ref{twUqY} states that the isomorphism~$\varphi$ of Theorem~\ref{mainthm} restricts to the twisted Yangian
$Y^{tw}(\mathfrak{g}_N)$ ($\mathfrak{g}_N = \mathfrak{o}_N$ or $\mathfrak{g}_N = \mathfrak{sp}_N$), viewed as
a~subalgebra of $Y(\mathfrak{gl}_N)$, and provides an isomorphism with a~certain subalgebra of
$\widetilde{Y}(\mathfrak{gl}_N)$ obtained from the generators of the corresponding twisted quantum loop
algebra.
Theorem~\ref{twUqY2} says that $Y^{tw}(\mathfrak{g}_N)$ is isomorphic to a~certain limit form of the twisted quantum
loop algebra and is analogous to Theorem~\ref{mainthm} for $Y(\mathfrak{gl}_N)$.
Results si\-mi\-lar to Theorems~\ref{twUqY} and~\ref{twUqY2} should hold for the twisted quantum loop algebras and twisted
Yangians of type AIII treated in~\cite{CGM} (resp.~\cite{MoRa}), and for the twisted quantum loop superalgebra and
twisted super-Yangians explored in~\cite{ChGu} (resp.~\cite{Na}).
When $\mathfrak{g}$ is a~classical simple Lie algebra other than $\mathfrak{s}\mathfrak{l}_n$, there is also a~way to
present its Yangian $Y(\mathfrak{g})$ using using an RTT-relation for some appropriate matrix~$R$, see~\cite{AMR}.
(In \textit{loc.
cit.}, $Y(\mathfrak{g})$ is obtained as a~quotient of what the authors call the extended Yangian $X(\mathfrak{g})$.) The
main results of~\cite{AMR} should admit analogues for quantum loop algebras, but the present authors ignore whether this
has been worked out or not.
If so, then Theorem~\ref{mainthm} below should have an analogue for Yangians and quantum loop algebras attached to
$\mathfrak{g}$ via RTT-relations when $\mathfrak{g}$ is of classical Dynkin type.
Furthermore, twisted Yangians which are coideal subalgebras of $Y(\mathfrak{g})$ and of $X(\mathfrak{g})$ have been
studied recently in~\cite{GuRe} and hopefully admit~$q$-analogues for which a~version of Theorem~\ref{twUqY2} should
hold also.
As for more general twisted Yangians associated to an arbitrary symmetric pair (see~\cite{Ma1} for instance), it is
reasonable to think that they are also limit forms of certain twisted quantum loop algebras: the quantum symmetric
Kac--Moody pairs introduced in~\cite{Ko} may be relevant here.

One motivation for this paper is that it should help in establishing an isomorphism using the RTT-generators between
completions of the Yangian and quantum loop algebra of $\mathfrak{gl}_n$ similar to the one obtained
in~\cite{GaTL1} using the generators from Drinfeld's second realization.
Once such a~result is obtained, one may hope that such an isomorphism would work also for twisted Yangians and twisted
quantum loop algebras.
Eventually, the goal would be to establish equiva\-lences between module categories for twisted Yangians and twisted
quantum loop algebras as in~\cite{GaTL2,GaTL3}.

\section[Quantum loop algebra and Yangian for $\mathfrak{gl}_N$]{Quantum loop algebra and Yangian for $\boldsymbol{\mathfrak{gl}_N}$}
%\label{glnsec}

We start by recalling the def\/inition of the Yangian and of the quantum loop algebra of $\mathfrak{gl}_N$ in
the RTT presentation~\cite{FRT}.
We view $\hbar$ as a~formal parameter.

\begin{definition}
\label{Ygln}
Let $R(u) = 1 \otimes 1 - \hbar u^{-1} \sum\limits_{i,j=1}^N E_{ij} \otimes E_{ji} $ where $E_{ij}$ denotes the usual
elementary matrix.
The Yangian $Y(\mathfrak{gl}_N)$ is the unital associative algebra over $\mathbb{C}[\hbar]$ generated~by
elements~$t_{ij}^{(r)}$ for~$r\in\mathbb{Z}_+$, $i,j \in \{1,\ldots,N\}$ satisfying $t_{ij}^{(0)}=\delta_{ij}$, together
with the following relations: if $t_{ij}(u)=\sum\limits_{r=0}^\infty {t_{ij}^{(r)}u^{-r}} \in
Y(\mathfrak{gl}_N)[[u^{-1}]]$ and $t(u) =\sum\limits_{i,j=1}^N {E_{ij} \otimes t_{ij}(u)}$, then
\begin{gather}
R(u-v)t_1(u)t_2(v)=t_2(v)t_1(u)R(u-v),
\label{RTT}
\end{gather}
where $t_1(u)$ (resp.\ $t_2(u)$) is the element of $\End_{\mathbb{C}}\big(\mathbb{C}^N\big) \otimes_{\mathbb{C}}
\End_{\mathbb{C}}\big(\mathbb{C}^N\big) \otimes_{\mathbb{C}} Y(\mathfrak{gl}_N)[[u^{-1}]]$ obtained by replacing
$E_{ij}$ by $E_{ij} \otimes I$ (resp.
$I \otimes E_{ij}$) in $t(u)$.
Here, we have also identif\/ied $R(u-v)$ with $R(u-v)\otimes 1$.
Relation~\eqref{RTT} is equivalent to the system of relations
\begin{gather*}
\big[t_{ij}^{(m+1)},t_{kl}^{(n)}\big]-\big[t_{ij}^{(m)},t_{kl}^{(n+1)}\big]
= \hbar\big(t_{kj}^{(m)}t_{il}^{(n)}-t_{kj}^{(n)}t_{il}^{(m)}\big).
\end{gather*}
\end{definition}

\begin{definition}
The quantum af\/f\/ine~$R$-matrix is the element of $\End_{\mathbb{C}}\big(\mathbb{C}^N\big) \otimes
\End_{\mathbb{C}}\big(\mathbb{C}^N\big) \otimes_{\mathbb{C}} \mathbb{C}[u,v]$ given~by
\begin{gather*}
R_q(u,v) = \sum\limits_{i,j=1}^N \big(uq^{-\delta_{ij}} - vq^{\delta_{ij}}\big) E_{ii} \otimes E_{jj} - \big(q-q^{-1}\big)u
\sum\limits_{\substack{i,j=1\\ i > j}}^N E_{ij} \otimes E_{ji}\\
\hphantom{R_q(u,v) =}{}
- \big(q-q^{-1}\big)v \sum\limits_{\substack{i,j=1\\ i < j}}^N E_{ij} \otimes E_{ji}.
\end{gather*}
\end{definition}
Set $\mathcal{L}(\mathfrak{gl}_N)=\mathfrak{gl}_N \otimes_{\mathbb{C}} \mathbb{C}[s,s^{-1}]$.
\begin{definition}
The quantum loop algebra $\mathfrak{U}_q(\mathcal{L}(\mathfrak{gl}_N))$ is the unital associative algebra
over~$\mathbb{C}(q)$ generated by $\{T_{ij}^{(r)},\overline{T}_{ij}^{(r)} \,|\, r\in\mathbb{Z}_+, \; i,j=1,\ldots,N\}$
which must satisfy the following relations
\begin{gather}
T_{ij}^{(0)} = 0 = \overline{T}_{ji}^{(0)} \qquad \text{if} \qquad  1 \leq i<j \leq N,
\nonumber
\\
T_{ii}^{(0)} \overline{T}_{ii}^{(0)} = 1 = \overline{T}_{ii}^{(0)} T_{ii}^{(0)} \qquad \forall\,   1 \leq i \leq N,
\nonumber
\\
R_q(u,v)T_1(u)T_2(v) = T_2(v) T_1(u) R_q(u,v),
\label{RTTql1}
\\
R_q(u,v)\overline{T}_1(u)\overline{T}_2(v) = \overline{T}_2(v) \overline{T}_1(u) R_q(u,v),
\label{RTTql2}
\\
R_q(u,v)\overline{T}_1(u) T_2(v) = T_2(v) \overline{T}_1(u) R_q(u,v),
\label{RTTql3}
\end{gather}
$T_a(u)$ and $\overline{T}_a(u)$ are obtained from $T(u):=\sum\limits_{i,j=1}^N{E_{ij} \otimes T_{ij}(u)}$ and
$\overline{T}(u):=\sum\limits_{i,j=1}^N{E_{ij} \otimes \overline{T}_{ij}(u)}$ in the same way as $t_a(u)$ in
Def\/inition~\ref{Ygln}; in this case, $T_{ij}(u)$, $\overline{T}_{ij}(u)$ are the elements
\begin{gather*}
T_{ij}(u)=\sum\limits_{r=0}^\infty{T_{ij}^{(r)}u^{-r}} \in
\mathfrak{U}_q(\mathcal{L}(\mathfrak{gl}_N))\big[\big[u^{-1}\big]\big],
\qquad
\overline{T}_{ij}(u)=\sum\limits_{r=0}^\infty{\overline{T}_{ij}^{(r)}u^{r}} \in
\mathfrak{U}_q(\mathcal{L}(\mathfrak{gl}_N))[[u]].
\end{gather*}
\end{definition}

Relations~\eqref{RTTql1},\eqref{RTTql2} and~\eqref{RTTql3} can be made more explicit in terms of the generators of
$\mathfrak{U}_q(\mathcal{L}(\mathfrak{gl}_N))$.
For instance,~\eqref{RTTql1} is equivalent to the following family of relations
\begin{gather}
\big(q^{-\delta_{ik}} T_{ij}^{(r+1)} T_{kl}^{(s)} - q^{\delta_{ik}}T_{ij}^{(r)} T_{kl}^{(s+1)}\big) - \big(q^{-\delta_{jl}}
T_{kl}^{(s)} T_{ij}^{(r+1)} - q^{\delta_{jl}} T_{kl}^{(s+1)} T_{ij}^{(r)}\big)
\nonumber
\\
\qquad{}
= \big(q-q^{-1}\big) \big(\delta_{i>k} T_{kj}^{(r+1)} T_{il}^{(s)} + \delta_{i<k} T_{kj}^{(r)} T_{il}^{(s+1)}\big)
\nonumber\\
\qquad\phantom{=}
- \big(q-q^{-1}\big)\big(\delta_{l>j} T_{kj}^{(s)} T_{il}^{(r+1)} + \delta_{l<j} T_{kj}^{(s+1)} T_{il}^{(r)}\big).
\label{RTTql12}
\end{gather}

Let $\mathcal{A}$ be the localization of $\mathbb{C}[q,q^{-1}]$ at the ideal $(q-1)$.
Let $\mathfrak{U}_{\mathcal{A}}(\mathcal{L}(\mathfrak{gl}_N))$ be the $\mathcal{A}$-subalgebra of
$\mathfrak{U}_q(\mathcal{L}(\mathfrak{gl}_N))$ generated by the elements $ \tau_{ij}^{(r)}$,
$\overline{\tau}_{ij}^{(r)}$ given~by
\begin{gather*}
\tau_{ij}^{(r)} = \frac{T_{ij}^{(r)}}{q-q^{-1}},
\qquad
\overline{\tau}_{ij}^{(r)} = \frac{\overline{T}_{ij}^{(r)}}{q-q^{-1}}
\qquad\text{for}\quad r\ge 0,\quad 1\le i,j \le N,
\end{gather*}
except that, when $r=0$ and $i=j$, we set
\begin{gather*}
\tau_{ii}^{(0)} = \frac{T_{ii}^{(0)}-1}{q-1}, \qquad \overline{\tau}_{ii}^{(0)} = \frac{\overline{T}_{ii}^{(0)}-1}{q-1}.
\end{gather*}

\begin{theorem}[Section~3 of~\cite{MRS}]
\label{qlim}
The assignment $E_{ij} s^r \mapsto \tau_{ij}^{(r)}$ $\forall\, r\ge 0$, $1\le i$, $j\le n$ except if $r=0$ and $1\le i<j\le n$,
$E_{ij} s^{-r} \mapsto -\overline{\tau}_{ij}^{(r)}$ $\forall\, r\ge 0$, $1\le i$, $j\le n$ except if $r=0$ and $1\le j<i\le n$,
induces an isomorphism $\mathfrak{U}(\mathcal{L}(\mathfrak{gl}_N)) \stackrel{\sim}{\longrightarrow}
\mathfrak{U}_{\mathcal{A}}(\mathcal{L}(\mathfrak{gl}_N)) \otimes_{\mathcal{A}} \mathbb{C}$, where $\mathbb{C}$
is viewed as an $\mathcal{A}$-module via $\mathcal{A}/(q-1) \stackrel{\sim}{\longrightarrow} \mathbb{C}$.
\end{theorem}

We have the following sequence of algebra homomorphisms similar to the one considered in~\cite{GuMa}
\begin{gather*}
\mathfrak{U}_{\mathcal{A}}(\mathcal{L}(\mathfrak{gl}_N)) \twoheadrightarrow
\mathfrak{U}_{\mathcal{A}}(\mathcal{L}(\mathfrak{gl}_N)) /
(q-1)\mathfrak{U}_{\mathcal{A}}(\mathcal{L}(\mathfrak{gl}_N)) \stackrel{\sim}{\longrightarrow}
\mathfrak{U}(\mathcal{L}(\mathfrak{gl}_N)) \stackrel{s\mapsto 1}{\twoheadrightarrow}
\mathfrak{U}(\mathfrak{gl}_N).
%\label{evalhom}
\end{gather*}

For $m\ge 0$, denote by $\mathsf{K}_m$ the Lie ideal of $\mathcal{L}(\mathfrak{gl}_N)$
spanned by $X\otimes s^r (s-1)^m \forall\, r\in\mathbb{Z}$, $\forall\, X\in\mathfrak{gl}_N$.
Let~$U$ be the subspace of $\mathfrak{U}_{\mathcal{A}}(\mathcal{L}(\mathfrak{gl}_N))$ spanned over
$\mathbb{C}$ by the generators $\tau_{ij}^{(r)}$, $\overline{\tau}_{ij}^{(r)}$ and observe that
 \begin{gather*}
 U \cap (q-1){\mathfrak{U}}_{\mathcal{A}}({\mathcal{L}}({\mathfrak{gl}}_N))
 = \operatorname{span}_{\mathbb{C}} \big\{\tau_{ii}^{(0)}
 +\overline{\tau}_{ii}^{(0)} \,\big|\, i=1,\ldots,N\big\}.
\end{gather*}
This is because $(q-1)\mathfrak{U}_{\mathcal{A}}(\mathcal{L}(\mathfrak{gl}_N))=\Ker(\psi)$,
where~$\psi$ is the composite
\begin{gather*}
\mathfrak{U}_{\mathcal{A}}(\mathcal{L}(\mathfrak{gl}_N)) \twoheadrightarrow
\mathfrak{U}_{\mathcal{A}}(\mathcal{L}(\mathfrak{gl}_N)) /
(q-1)\mathfrak{U}_{\mathcal{A}}(\mathcal{L}(\mathfrak{gl}_N)) \stackrel{\sim}{\longrightarrow}
\mathfrak{U}(\mathcal{L}(\mathfrak{gl}_N)).
\end{gather*}

Moreover, notice that
\begin{gather*}
\tau_{ii}^{(0)}+\overline{\tau}_{ii}^{(0)} = \frac{\overline{T}_{ii}^{(0)}\big(T_{ii}^{(0)}-1\big)^2}{q-1} = (q-1)
\big((q-1)\overline{\tau}_{ii}^{(0)}+1\big) \big(\tau_{ii}^{(0)}\big)^2 \in
(q-1)\mathfrak{U}_{\mathcal{A}}(\mathcal{L}(\mathfrak{gl}_N)).
\end{gather*}
Let $\mathbf{K}_0 = \mathfrak{U}_{\mathcal{A}}(\mathcal{L}(\mathfrak{gl}_N))$.
For $m \geq 1$, let $\mathbb{K}_m$ be the two-sided ideal of
$\mathfrak{U}_{\mathcal{A}}(\mathcal{L}(\mathfrak{gl}_N))$ generated by $\psi^{-1}(\mathsf{K}_m) \cap U$, and
set $\mathbf{K}_m$ to be the sum of the ideals $(q-q^{-1})^{m_0}\mathbb{K}_{m_1} \cdots \mathbb{K}_{m_k}$ with $m_0 + m_1
+ \cdots + m_k \ge m$.
This is slightly dif\/ferent from the def\/inition of the analogous ideals $\mathbf{K}_m$ in~\cite{GuMa} in the case of
$\mathfrak{s}\mathfrak{l}_N$ because, for the $\mathfrak{gl}_N$ case, $\mathbf{K}_1^m$ is strictly smaller
than $\mathbf{K}_m$.

Let $\widetilde{Y}(\mathfrak{gl}_N)$ be the $\mathbb{C}$-algebra
\begin{gather*}
\widetilde{Y}(\mathfrak{gl}_N) = \bigoplus_{m=0}^{\infty} \mathbf{K}_{m}/\mathbf{K}_{m+1}.
\end{gather*}
$\widetilde{Y}(\mathfrak{gl}_N)$ can be viewed as a~$\mathbb{C}[\hbar]$-algebra if we set
$\hbar=\overline{q-q^{-1}} \in \mathbf{K}_1/\mathbf{K}_2$.
Note that the f\/irst quotient is isomorphic to $\mathfrak{U}(\mathfrak{gl}_N)$ by def\/inition.

\begin{theorem}
\label{mainthm}
$ \widetilde{Y}(\mathfrak{gl}_N)$ is isomorphic to $Y(\mathfrak{gl}_N)$.
\end{theorem}

For the analogue of our f\/irst theorem for an arbitrary complex semisimple Lie algebra, see~\cite{Dr1} and~\cite{GuMa}.

For $m,r\ge 0$, def\/ine recursively elements $T_{ij}^{(r,m)}$ in the following way
\begin{gather*}
T_{ij}^{(r,0)} = \tau_{ij}^{(r)} \qquad\text{and}\qquad T_{ij}^{(r,m+1)} = T_{ij}^{(r+1,m)} - T_{ij}^{(r,m)},
\end{gather*}
except that if $i<j$, $T_{ij}^{(0,0)} = -\overline{\tau}_{ij}^{(0)}$.

{\sloppy One can check by induction on~$m$ that, for every~$r$, $\psi\big(T_{ij}^{(r,m)}\big)=E_{ij}s^r(s-1)^m$, hence $T_{ij}^{(r,m)}
\in \mathbf{K}_m$.
Set $\xi_{ij}^{(r,m)} = \overline{T_{ij}^{(r,m)}} \in \mathbf{K}_m/\mathbf{K}_{m+1}$.

}

\begin{proof}[Proof of Theorem~\ref{mainthm}] We will prove that there is an isomorphism $\varphi\colon Y(\mathfrak{gl}_N)
\stackrel{\sim}{\longrightarrow} \widetilde{Y}(\mathfrak{gl}_N)$ given by $t_{ij}^{(m+1)} \mapsto
\xi_{ij}^{(0,m)}$ for $m\ge 0$.

Relation~\eqref{RTTql12} can be rewritten in the following way
\begin{gather*}
q^{-\delta_{ik}}\big(T_{ij}^{(r+1)} - T_{ij}^{(r)}\big) T_{kl}^{(s)}- T_{ij}^{(r)} \big(q^{\delta_{ik}}T_{kl}^{(s+1)}
- q^{-\delta_{ik}}T_{kl}^{(s)}\big)
\\
\qquad\phantom{=}{}
- \big(q^{-\delta_{jl}}T_{kl}^{(s)} \big(T_{ij}^{(r+1)} - T_{ij}^{(r)}\big) - \big(q^{\delta_{jl}}T_{kl}^{(s+1)} -
q^{-\delta_{jl}}T_{kl}^{(s)}\big) T_{ij}^{(r)} \big)
\\
\qquad
= \big(q-q^{-1}\big) \big(\delta_{i>k}T_{kj}^{(r+1)} T_{il}^{(s)} + \delta_{i<k} T_{kj}^{(r)} T_{il}^{(s+1)}\big)
\\
\qquad\phantom{=}{}
  - \big(q-q^{-1}\big) \big(\delta_{l>j} T_{kj}^{(s)} T_{il}^{(r+1)} + \delta_{l<j} T_{kj}^{(s+1)} T_{il}^{(r)} \big).
\end{gather*}

If $r,s \geq 1$ then, after rearranging and dividing both sides by $(q-q^{-1})^2$, we get
\begin{gather*}
q^{-\delta_{ik}} \big(T_{ij}^{(r,1)} T_{kl}^{(s,0)}- T_{ij}^{(r,0)} T_{kl}^{(s,1)} \big) - \big(q^{\delta_{ik}} -
q^{-\delta_{ik}}\big) T_{ij}^{(r,0)} T_{kl}^{(s+1,0)}
\\
\qquad\phantom{=}{}
 - q^{-\delta_{jl}} \big(T_{kl}^{(s,0)} T_{ij}^{(r,1)} - T_{kl}^{(s,1)} T_{ij}^{(r,0)} \big) + \big(q^{\delta_{jl}} -
q^{-\delta_{jl}}\big) T_{kl}^{(s+1,0)} T_{ij}^{(r,0)}
\\
\qquad
= \big(q-q^{-1}\big) \big(  \delta_{i>k} T_{kj}^{(r+1,0)} T_{il}^{(s,0)} + \delta_{i<k} T_{kj}^{(r,0)} T_{il}^{(s+1,0)}\big)
\\
\qquad\phantom{=}{}
  - \big(q-q^{-1}\big) \big(\delta_{l>j} T_{kj}^{(s,0)} T_{il}^{(r+1,0)} + \delta_{l<j} T_{kj}^{(s+1,0)} T_{il}^{(r,0)} \big).
\end{gather*}

Using $T_{ij}^{(r+1,m)} - T_{ij}^{(r,m)}=T_{ij}^{(r,m+1)}$ and $T_{kl}^{(s+1,n)} - T_{kl}^{(s,n)}=T_{kl}^{(s,n+1)}$, we
deduce by induction on~$m$ and~$n$ that, for all $r,s \geq 1$ and all $m,n\ge 0$,
\begin{gather}
q^{-\delta_{ik}} \big(T_{ij}^{(r,m+1)}T_{kl}^{(s,n)} - T_{ij}^{(r,m)} T_{kl}^{(s,n+1)} \big) - \big(q^{\delta_{ik}} -
q^{-\delta_{ik}}\big) T_{ij}^{(r,m)} T_{kl}^{(s+1,n)}
\nonumber
\\
\qquad\phantom{=}{}
  - q^{-\delta_{jl}} \big(T_{kl}^{(s,n)} T_{ij}^{(r,m+1)} - T_{kl}^{(s,n+1)} T_{ij}^{(r,m)} \big) + \big(q^{\delta_{jl}} -
q^{-\delta_{jl}}\big) T_{kl}^{(s+1,n)} T_{ij}^{(r,m)}
\nonumber
\\
\qquad
= \big(q-q^{-1}\big) \big(  \delta_{i>k} T_{kj}^{(r+1,m)} T_{il}^{(s,n)} + \delta_{i<k} T_{kj}^{(r,m)} T_{il}^{(s+1,n)}\big)
\nonumber
\\
\qquad\phantom{=}{}
  - \big(q-q^{-1}\big) \big(\delta_{l>j} T_{kj}^{(s,n)} T_{il}^{(r+1,m)} + \delta_{l<j} T_{kj}^{(s+1,n)} T_{il}^{(r,m)} \big).
\label{rs}
\end{gather}
Consider the case $r=s=1$ in~\eqref{rs}.
Using $T_{ij}^{(r+1,m)}=T_{ij}^{(r,m+1)}+T_{ij}^{(r,m)}$ and $T_{kl}^{(r+1,n)}=T_{kl}^{(r,n+1)}+T_{kl}^{(r,n)}$ for all
$r\ge 0$, we obtain, for all $m,n\ge 0$,
\begin{gather*}
  q^{-\delta_{ik}} \big(T_{ij}^{(0,m+2)} + T_{ij}^{(0,m+1)}\big) T_{kl}^{(0,n+1)} -
q^{-\delta_{ik}}\big(T_{ij}^{(0,m+1)} + T_{ij}^{(0,m)}\big) T_{kl}^{(0,n+2)}
\\
\qquad\phantom{=}{}
  - \big(q^{\delta_{ik}} - q^{-\delta_{ik}}\big) \big(T_{ij}^{(0,m+1)} + T_{ij}^{(0,m)}\big) T_{kl}^{(1,n+1)} + q^{-\delta_{ik}}
\big(T_{ij}^{(0,m+2)} + T_{ij}^{(0,m+1)}\big) T_{kl}^{(0,n)}
\\
\qquad\phantom{=}{}
  - q^{-\delta_{ik}}\big(T_{ij}^{(0,m+1)} + T_{ij}^{(0,m)}\big) T_{kl}^{(0,n+1)} - \big(q^{\delta_{ik}} -
q^{-\delta_{ik}}\big)\big(T_{ij}^{(0,m+1)} + T_{ij}^{(0,m)}\big) T_{kl}^{(1,n)}
\\
\qquad\phantom{=}{}
  - q^{-\delta_{jl}} T_{kl}^{(0,n+1)} \big(T_{ij}^{(0,m+2)} + T_{ij}^{(0,m+1)}\big) +
q^{-\delta_{jl}}T_{kl}^{(0,n+2)}\big(T_{ij}^{(0,m+1)} + T_{ij}^{(0,m)}\big)
\\
\qquad\phantom{=}{}
  + \big(q^{\delta_{jl}} - q^{-\delta_{jl}}\big) T_{kl}^{(1,n+1)} \big(T_{ij}^{(0,m+1)} + T_{ij}^{(0,m)}\big) - q^{-\delta_{jl}}
T_{kl}^{(0,n)} \big(T_{ij}^{(0,m+2)} + T_{ij}^{(0,m+1)}\big)
\\
\qquad\phantom{=}{}
  + q^{-\delta_{jl}}T_{kl}^{(0,n+1)} \big(T_{ij}^{(0,m+1)} + T_{ij}^{(0,m)}\big) + \big(q^{\delta_{jl}} - q^{-\delta_{jl}}\big)
T_{kl}^{(1,n)} \big(T_{ij}^{(0,m+1)} + T_{ij}^{(0,m)}\big)
\\
\qquad
  = \big(q-q^{-1}\big) \Big(\delta_{i>k} \big(T_{kj}^{(1,m+1)} + T_{kj}^{(1,m)}\big) T_{il}^{(0,n+1)} + \delta_{i<k}
\big(T_{kj}^{(0,m+1)} + T_{kj}^{(0,m)}\big) T_{il}^{(1,n+1)} \Big)
\\
\qquad\phantom{=}{}
+ \big(q-q^{-1}\big) \Big(\delta_{i>k} \big(T_{kj}^{(1,m+1)} + T_{kj}^{(1,m)}\big) T_{il}^{(0,n)} + \delta_{i<k}
\big(T_{kj}^{(0,m+1)} + T_{kj}^{(0,m)}\big) T_{il}^{(1,n)}\Big)
\\
\qquad\phantom{=}{}
- \big(q-q^{-1}\big) \Big(\delta_{l>j} T_{kj}^{(0,n+1)} \big(T_{il}^{(1,m+1)} + T_{il}^{(1,m)}\big) + \delta_{l<j}
T_{kj}^{(1,n+1)} \big(T_{il}^{(0,m+1)} + T_{il}^{(0,m)}\big) \Big)
\\
\qquad\phantom{=}{}
- \big(q-q^{-1}\big) \Big(\delta_{l>j} T_{kj}^{(0,n)} \big(T_{il}^{(1,m+1)} + T_{il}^{(1,m)}\big) + \delta_{l<j} T_{kj}^{(1,n)}
\big(T_{il}^{(0,m+1)} + T_{il}^{(0,m)}\big) \Big).
\end{gather*}

Notice that both sides of this last equality are in $\mathbf{K}_{m+n+1}$ (and some of the terms are even in
$\mathbf{K}_{m+n+2}$ or in $\mathbf{K}_{m+n+3}$).
Modulo $\mathbf{K}_{m+n+2}$, we obtain the congruence
\begin{gather*}
q^{-\delta_{ik}} \big(T_{ij}^{(0,m+1)} T_{kl}^{(0,n)}- T_{ij}^{(0,m)} T_{kl}^{(0,n+1)} \big) - \big(q^{\delta_{ik}} -
q^{-\delta_{ik}}\big) T_{ij}^{(0,m)} T_{kl}^{(1,n)}
\\
\qquad\phantom{=}{}
  - q^{-\delta_{jl}} \big(T_{kl}^{(0,n)} T_{ij}^{(0,m+1)} - T_{kl}^{(0,n+1)} T_{ij}^{(0,m)} \big) + \big(q^{\delta_{jl}} -
q^{-\delta_{jl}}\big) T_{kl}^{(1,n)} T_{ij}^{(0,m)}
\\
\qquad
\equiv \big(q-q^{-1}\big) \big(  \delta_{i>k} T_{kj}^{(1,m)} T_{il}^{(0,n)} + \delta_{i<k} T_{kj}^{(0,m)} T_{il}^{(1,n)}\big)
\\
\qquad\phantom{=}{}
  - \big(q-q^{-1}\big) \big(\delta_{l>j} T_{kj}^{(0,n)} T_{il}^{(1,m)} + \delta_{l<j} T_{kj}^{(1,n)} T_{il}^{(0,m)}\big).
\end{gather*}
Moreover, modulo $\mathbf{K}_{m+n+2}$, we also have
\begin{gather*}
\big(q^{\delta_{ik}} - q^{-\delta_{ik}}\big) T_{ij}^{(0,m)} T_{kl}^{(1,n)} \equiv \big(q^{\delta_{ik}} - q^{-\delta_{ik}}\big)
T_{ij}^{(0,m)} T_{kl}^{(0,n)},
\\
\big(q^{\delta_{jl}}-q^{-\delta_{jl}}\big) T_{kl}^{(1,n)} T_{ij}^{(0,m)} \equiv \big(q^{\delta_{jl}}-q^{-\delta_{jl}}\big)
T_{kl}^{(0,n)} T_{ij}^{(0,m)},
\\
\big(q-q^{-1}\big) T_{kj}^{(1,m)} T_{il}^{(0,n)} \equiv \big(q-q^{-1}\big)T_{kj}^{(0,m)} T_{il}^{(0,n)},
\\
\big(q-q^{-1}\big) T_{kj}^{(0,m)}
T_{il}^{(1,n)} \equiv \big(q-q^{-1}\big) T_{kj}^{(0,m)} T_{il}^{(0,n)},
\\
\big(q-q^{-1}\big) T_{kj}^{(0,n)} T_{il}^{(1,m)} \equiv \big(q-q^{-1}\big) T_{kj}^{(0,n)} T_{il}^{(0,m)},
\\
\big(q-q^{-1}\big) T_{kj}^{(1,n)}
T_{il}^{(0,m)} \equiv \big(q-q^{-1}\big) T_{kj}^{(0,n)} T_{il}^{(0,m)}.
\end{gather*}
Therefore, passing to the quotient $\mathbf{K}_{m+n+1}/\mathbf{K}_{m+n+2}$, we obtain
\begin{gather*}
\big(\xi_{ij}^{(0,m+1)} \xi_{kl}^{(0,n)} - \xi_{ij}^{(0,m)} \xi_{kl}^{(0,n+1)} \big)- \delta_{ik} \hbar
\xi_{ij}^{(0,m)} \xi_{kl}^{(0,n)}
\\
\qquad\phantom{=}{}
  - \big(\xi_{kl}^{(0,n)} \xi_{ij}^{(0,m+1)} - \xi_{kl}^{(0,n+1)} \xi_{ij}^{(0,m)} \big) + \delta_{jl} \hbar
\xi_{kl}^{(0,n)} \xi_{ij}^{(0,m)}
\\
\qquad
= \hbar \big(\delta_{i>k} \xi_{kj}^{(0,m)} \xi_{il}^{(0,n)}+ \delta_{i<k} \xi_{kj}^{(0,m)} \xi_{il}^{(0,n)}\big)
  - \hbar \big(\delta_{l>j} \xi_{kj}^{(0,n)} \xi_{il}^{(0,m)} + \delta_{l<j} \xi_{kj}^{(0,n)} \xi_{il}^{(0,m)}\big).
\end{gather*}
This last relation is equivalent to
\begin{gather*}
\big[\xi_{ij}^{(0,m+1)}, \xi_{kl}^{(0,n)}\big] - \big[\xi_{ij}^{(0,m)}, \xi_{kl}^{(0,n+1)}\big] = \hbar \big(\xi_{kj}^{(0,m)}
\xi_{il}^{(0,n)} - \xi_{kj}^{(0,n)} \xi_{il}^{(0,m)}\big).
\end{gather*}
This holds for all $m,n\ge 0$.

All the previous computations prove that $\varphi\colon Y(\mathfrak{gl}_N) \longrightarrow
\widetilde{Y}(\mathfrak{gl}_N)$ given by $\varphi(t_{ij}^{(m+1)}) = \xi_{ij}^{(0,m)}$ for $m\ge 0$ is an
algebra homomorphism.
We still have to show that~$\varphi$ is injective and surjective.

We will f\/irst demonstrate surjectivity.
Towards this end, we def\/ine elements $\overline{T}_{ij}^{(r,m)}$ as follows.
Let $\overline{T}_{ij}^{(r,0)}=\overline{\tau}_{ij}^{(r)}$, except that $\overline{T}_{ij}^{(0,0)}=-\tau_{ij}^{(0)}$
when $i \geq j$.
Then, for each $m \geq 0$, let
\begin{gather*}
\overline{T}_{ij}^{(r,m+1)}=\overline{T}_{ij}^{(r+1,m)}-\overline{T}_{ij}^{(r,m)}.
\end{gather*}

Also for each $m \geq 0$, let $\widetilde{T}_{ij}^{(0,m)}=\overline{T}_{ij}^{(0,m)}$ and
$\widetilde{T}_{ij}^{(m,m)}=(-1)^{m+1}T_{ij}^{(0,m)}$, and for $1 \leq r \leq m$ def\/ine recursively
\begin{gather*}
\widetilde{T}_{ij}^{(r,m+1)}=\widetilde{T}_{ij}^{(r-1,m)}-\widetilde{T}_{ij}^{(r,m)}.
\end{gather*}

Induction on~$m$ shows that the elements $T_{ij}^{(r,m)}$, $\overline{T}_{ij}^{(r,m)}$ and $\widetilde{T}_{ij}^{(r,m)}$
respectively map via~$\psi$ to the elements $E_{ij}s^r(s-1)^m$, $(-1)^{m+1}E_{ij}s^{-(m+r)}(s-1)^m$ and
$(-1)^{m+1}E_{ij}s^{-(m-r)}(s-1)^m$ in $\mathfrak{U}(\mathcal{L}(\mathfrak{gl}_N))$.
It follows that, for f\/ixed~$m$, the images of those elements under~$\psi$ span $\mathsf{K}_m$.
Moreover, all those elements are in~$U$ by def\/inition.

Note that for any f\/ixed $X \in \psi^{-1}(\mathsf{K}_m) \cap U$, there exists some element~$Y$ in
 \begin{gather*}
\operatorname{span}_{\mathbb{C}}\big\{T_{ij}^{(r,m)},\overline{T}_{ij}^{(r,m)},\widetilde{T}_{ij}^{(r,m)} \,|\, i,j=1,\ldots,N,\; r \in
\mathbb{Z}_+\big\}
\end{gather*}
for which $X-Y \in (q-1)\mathfrak{U}_{\mathcal{A}}(\mathcal{L}(\mathfrak{gl}_N))$, because the map
$\mathfrak{U}_{\mathcal{A}}\mathcal{L}(\mathfrak{gl}_N)
/(q-1)\mathfrak{U}_{\mathcal{A}}(\mathcal{L}(\mathfrak{gl}_N))
\stackrel{\sim}{\longrightarrow} \mathfrak{U}(\mathcal{L}(\mathfrak{gl}_N))$ is an isomorphism.
Since $X-Y$ is also in~$U$, we have{\samepage
\begin{gather*}
X-Y \in \operatorname{span}_{\mathbb{C}} \big\{\tau_{ii}^{(0)}+\overline{\tau}_{ii}^{(0)} \,|\, i=1,\ldots,N\big\} \subset \mathbf{K}_{\ell}
 \; \forall\, \ell \geq 0.
\end{gather*}
That $\tau_{ii}^{(0)}+\overline{\tau}_{ii}^{(0)}$ is in $\mathbf{K}_{\ell}$ for all $\ell \geq 0$ is a~consequence of the
fact that $\psi(\tau_{ii}^{(0)}+\overline{\tau}_{ii}^{(0)})=0$.}

It follows that any element of $\mathbf{K}_m$ is congruent modulo $\mathbf{K}_{m+1}$ to a~sum of monomials of the form
$f(q)(q-q^{-1})^{m_0}\mathcal{M}$ where $f(q) \in \mathcal{A}$ is not divisible by $q-1$ and $\mathcal{M}=\tau_{i_1
j_1}^{(r_1,m_1)} \cdots \tau_{i_k j_k}^{(r_k,m_k)}$ with
\begin{gather*}
\tau_{i_d j_d}^{(r_d,m_d)} \in \big\{T_{i_d j_d}^{(r_d,m_d)}, \overline{T}_{i_d j_d}^{(r_d,m_d)}, \widetilde{T}_{i_d
j_d}^{(r_d,m_d)}\big\}
\end{gather*}
and $m_0+\cdots+m_k \geq m$.
Moreover, since $T_{i_d j_d}^{(r_d,m_d)}-T_{i_d j_d}^{(r_d-1,m_d)}=T_{i_d j_d}^{(r_d-1,m_d+1)} \in \mathbb{K}_{m_d+1}$
(and similarly for $\overline{T}_{i_d j_d}^{(r_d,m_d)}$ and $\widetilde{T}_{i_d j_d}^{(r_d,m_d)}$), we can reduce modulo
$\mathbf{K}_{m+1}$ to the case when $r_d=0$ for each $d=1,\ldots,k$.

Observe that modulo $\mathbb{K}_{m_d+1}$, we have
\begin{gather}
\overline{T}_{ij}^{(0,m_d)} \equiv (-1)^{m_d+1}T_{ij}^{(0,m_d)} \equiv \widetilde{T}_{ij}^{(0,m_d)}.
\label{TbarT}
\end{gather}
To see this, just take the dif\/ference of the elements on each side (this dif\/ference is in~$U$ by def\/inition) and
apply~$\psi$.
Finally, observe that we can replace $f(q)$ by $f(1)$ modulo $\mathbf{K}_{1}$.

In summary, we have shown that each of the monomials $f(q)(q-q^{-1})^{m_0}\mathcal{M}$ in $\mathbf{K}_m$ is congruent
modulo $\mathbf{K}_{m+1}$ to
\begin{gather*}
f(1)\big(q-q^{-1}\big)^{m_0} T_{i_1 j_1}^{(0,m_1)} \cdots T_{i_k,j_k}^{(0,m_k)}
\end{gather*}
up to a~sign.
The image modulo $\mathbf{K}_{m+1}$ of this element is
\begin{gather*}
f(1)\hbar^{m_0}\xi_{i_1 j_1}^{(0,m_1)} \cdots \xi_{i_k j_k}^{(0,m_k)}
\end{gather*}
and this is in the image of~$\varphi$ by def\/inition.
This completes the proof that~$\varphi$ is surjective.

To prove that~$\varphi$ is injective, it is enough to show that the basis of $Y(\mathfrak{gl}_N)$ given~by
ordered monomials (for some f\/ixed order) in the generators $t_{ij}^{(m)}$ is mapped via~$\varphi$ to some linearly
independent set in $\widetilde{Y}(\mathfrak{gl}_N)$.

By def\/inition, any two ordered monomials $\hbar^{m_0}t_{i_1 j_1}^{(m_1+1)} \cdots t_{i_a j_a}^{(m_a+1)}$ and
$\hbar^{n_0}t_{k_1 l_1}^{(n_1+1)} \cdots t_{k_b l_b}^{(n_b+1)}$ with each $m_d, n_d \geq 0$ and $m_0+m_1+\cdots+m_a \neq
n_0+n_1 + \cdots + n_b$ are mapped via~$\varphi$ to distinct graded pieces in $\widetilde{Y}(\mathfrak{gl}_N)$.
It therefore suf\/f\/ices to show that for each f\/ixed~$m$, the images under~$\varphi$ of all the ordered monomials
$\hbar^{m_0} t_{i_1 j_1}^{(m_1+1)} \cdots t_{i_a j_a}^{(m_a+1)}$ with $m_1+\cdots+m_a=m$ are linearly independent in
$\mathbf{K}_m/\mathbf{K}_{m+1}$.

Consider any linear combination over $\mathbb{C}$ of ordered monomials $(\overline{q-q^{-1}})^{m_0}\xi_{i_1
j_1}^{(0,m_1)} \cdots \xi_{i_a j_a}^{(0,m_a)}$ with $m_0+m_1+\cdots+m_a=m$, and suppose that this sum is zero in
$\mathbf{K}_m/\mathbf{K}_{m+1}$.
Then we have a~linear combination~$S$ of ordered monomials $(q-q^{-1})^{m_0}T_{i_1 j_1}^{(0,m_1)} \cdots T_{i_a
j_a}^{(0,m_a)}$ which is not just in~$\mathbf{K}_m$, but also in~$\mathbf{K}_{m+1}$.
We can assume that the minimum of the exponents~$m_0$ is~$0$.

$\psi(S)$ is a~linear combination of some of the ordered monomials $E_{i_1 j_1}(s-1)^{m_1} \cdots E_{i_a j_a}
(s-1)^{m_a}$.
On the other hand, since~$S$ is also in $\mathbf{K}_{m+1}$, $\psi(S)$ can be expressed as a~linear combination of
monomials of the form $E_{k_1 l_1} s^{r_1} (s-1)^{n_1} \cdots E_{k_b l_b} s^{r_b} (s-1)^{n_b}$ with $r_1,\ldots,r_b \in
\mathbb{Z}$ and $n_1+\cdots+n_b \geq m+1$.
If $r_1 = \cdots = r_b=0$, then this is impossible unless the coef\/f\/icients of both linear combinations vanish.
Let's prove that the same is true more generally.

For each $r \geq 1$, we have a~composite of algebra homomorphisms
\begin{gather*}
\mathfrak{U}(\mathcal{L}(\mathfrak{gl}_N)) \stackrel{\Delta}{\to}
\mathfrak{U}(\mathcal{L}(\mathfrak{gl}_N))^{\otimes r} \stackrel{f^{\otimes r}}{\rightarrow}
\End_{\mathbb{C}}\big(\mathbb{C}^N\big)^{\otimes r} \otimes \mathbb{C}\big[x_1^{\pm 1},x_2^{\pm 1},\ldots,x_r^{\pm 1}\big],
\end{gather*}
where~$\Delta$ is the standard coproduct on the enveloping algebra of a~Lie algebra and the homomorphism $f\colon
\mathfrak{U}(\mathcal{L}(\mathfrak{gl}_N)) \to \End_{\mathbb{C}}\big(\mathbb{C}^N\big) \otimes
\mathbb{C}[x,x^{-1}]$ is given by $f(E_{ij}x^t) = E_{ij} \otimes x^t$.

We also have for each choice of~$r$ nonnegative integers $\alpha_1, \ldots, \alpha_r$ a~dif\/ferential operator
\begin{gather*}
\partial_{\alpha_1,\ldots,\alpha_r}\colon \ \End_{\mathbb{C}}\big(\mathbb{C}^N\big)^{\otimes r} \otimes \mathbb{C}\big[x_1^{\pm
1},x_2^{\pm 1},\ldots,x_r^{\pm 1}\big] \to \End_{\mathbb{C}}\big(\mathbb{C}^N\big)^{\otimes r}
\end{gather*}
given~by
\begin{gather*}
\partial_{\alpha_1, \ldots, \alpha_r}=\left.\frac{\partial^{\alpha_1}}{\partial x_1^{\alpha_1}}
\frac{\partial^{\alpha_2}}{\partial x_2^{\alpha_2}} \cdots \frac{\partial^{\alpha_r}}{\partial x_r^{\alpha_r}}
\right\rvert_{x_1,\ldots,x_r=1}.
\end{gather*}
Take $r \geq \max\{a,b\}$ where the maximum is taken over all the monomials in $\psi(S)$ (with~$a$ and~$b$ related to
the monomials in~$S$ and $\psi(S)$ as above) and note that for any choice of $\alpha_1,\ldots,\alpha_r$ with $\alpha_1+
\cdots + \alpha_r=m$, $\psi(S)$ is in the kernel of the composite $\partial_{\alpha_1,\ldots,\alpha_r} \circ f^{\otimes
r} \circ \Delta$ because $n_1 + \cdots + n_b \geq m+1$.

On the other hand, if~$S$ is nonzero, then we can f\/ind some $\alpha_1,\ldots,\alpha_r$ such that $\alpha_1 + \cdots +
\alpha_r=m$ and $\psi(S)$ is not in the kernel of $\partial_{\alpha_1,\ldots,\alpha_r} \circ f^{\otimes r} \circ
\Delta$: just choose any of the ordered monomials $E_{i_1 j_1}(s-1)^{m_1} \cdots E_{i_a j_a} (s-1)^{m_a}$ in $\psi(S)$
and set $\alpha_1=m_1, \ldots, \alpha_a=m_a$ and $\alpha_d=0$ for $d>a$.

We have just obtained a~contradiction, so $S=0$ and the linear sum of ordered monomials $\xi_{i_1 j_1}^{(0,m_1)} \cdots
\xi_{i_a j_a}^{(0,m_a)}$ must in fact be trivial, as desired.
\end{proof}

\section[Twisted quantum loop algebras and Yangians of type $\mathfrak{o}_N$ and $\mathfrak{sp}_N$]{Twisted quantum loop algebras and Yangians\\ of type $\boldsymbol{\mathfrak{o}_N}$ and $\boldsymbol{\mathfrak{sp}_N}$}
\label{twsec}
We will now prove an analogue of Theorem~\ref{mainthm} for certain twisted Yangians and twisted quantum loop algebras
associated to the symmetric pairs $(\mathfrak{gl}_N,\mathfrak{o}_N)$ and
$(\mathfrak{gl}_N,\mathfrak{sp}_N)$.
We will treat these two cases simultaneously and denote by $\mathfrak{g}_N$ either $\mathfrak{o}_N$ or
$\mathfrak{sp}_N$.
Whenever we use the symbols $\pm$ or $\mp$, it is understood that the sign on the top is used for
$\mathfrak{g}_N=\mathfrak{o}_N$ and the sign on the bottom is used for $\mathfrak{g}_N=\mathfrak{sp}_N$.
In this section, we will use~$t$ to denote transposition in the f\/irst factor of a~tensor product of matrices.

In the orthogonal case, let $G=(g_{ij})$ be the $N \times N$ identity matrix.
For the symplectic case, we take
\begin{gather*}
G=\sum\limits_{k=1}^{N/2} (E_{2k-1,2k}-E_{2k,2k-1}),
\end{gather*}
which makes sense since $\mathfrak{sp}_N$ is only def\/ined when~$N$ is even.
Similarly, let $B=(b_{ij})$ be the $N \times N$ identity matrix in the orthogonal case, while in the symplectic case we
take
\begin{gather*}
B=\sum\limits_{k=1}^{N/2} (qE_{2k-1,2k}-E_{2k,2k-1}).
\end{gather*}

\begin{definition}
\label{Ytw}
Let $R(u)$ be as given in Def\/inition~\ref{Ygln}.
The twisted Yangian $Y^{tw}(\mathfrak{g}_N)$ is the unital associative algebra over $\mathbb{C}[\hbar]$ generated~by
$\big\{s_{ij}^{(r)} \,|\, r\in\mathbb{Z}_+,\; i,j=1,\ldots,N\big\}$ where $s_{ij}^{(0)}=g_{ij}$, together with the following
relations: if $s_{ij}(u)=\sum\limits_{r=0}^\infty s_{ij}^{(r)}u^{-r}$ and $s(u)=\sum\limits_{i,j=1}^N E_{ij} \otimes
s_{ij}(u)$, then
\begin{gather*}
R(u-v)s_1(u)R^t(-u-v)s_2(v)=s_2(v)R^t(-u-v)s_1(u)R(u-v)
\end{gather*}
and
\begin{gather*}
s^t(-u)= \pm s(u) + \hbar\frac{s(u)-s(-u)}{2u},
\end{gather*}
where $s_1(u)$ (resp.\
$s_2(u)$) is the element of $\End_{\mathbb{C}}\big(\mathbb{C}^N\big) \otimes \End_{\mathbb{C}}\big(\mathbb{C}^N\big) \otimes
Y^{tw}(\mathfrak{g}_N)[[u^{-1}]]$ obtained by replacing $E_{ij}$ by $E_{ij} \otimes I$ (resp.
$I \otimes E_{ij}$) in $s(u)$.
\end{definition}

The twisted Yangian is a~deformation of the universal enveloping algebra of the twisted current algebra
$\mathfrak{g}_N^{tw}[s]$ which is def\/ined in the following way.

\begin{definition}
%\label{twcurrent}
Let~$\sigma$ be the automorphism of $\mathfrak{gl}_N$ given~by
\begin{gather}
\label{orthauto}
\sigma(E_{ij})=-E_{ji}
\end{gather}
if $\mathfrak{g}_N=\mathfrak{o}_N$, while
\begin{gather}
\label{sympauto}
\sigma(E_{ij})=(-1)^{i+j-1}E_{j'i'}
\end{gather}
if $\mathfrak{g}_N=\mathfrak{sp}_N$.
Here, $i'=i-1$ if~$i$ is even and $i'=i+1$ is~$i$ is odd.
The twisted current algebra is the subalgebra of $\mathfrak{gl}_N[s]$ given~by
\begin{gather*}
\mathfrak{g}_N^{tw}[s] = \{A(s) \in \mathfrak{gl}_N[s] \,|\, \sigma(A(s))=A(-s)\}.
\end{gather*}
\end{definition}

The twisted Yangians can be regarded as subalgebras of the Yangian for $\mathfrak{gl}_N$:

\begin{proposition}[\cite{Ol}]
\label{embY}
The assignment
\begin{gather*}
s_{ij}^{(r)} \mapsto \sum\limits_{k=1}^N \left(g_{kj}t_{ik}^{(r)} + (-1)^r g_{ik} t_{jk}^{(r)} \right) + \hbar
\sum\limits_{k,l=1}^N \sum\limits_{p=1}^{r-1} (-1)^{r-p} g_{kl} t_{ik}^{(p)} t_{jl}^{(r-p)}
\end{gather*}
provides an embedding of $Y^{tw}(\mathfrak{g}_N)$ into $Y(\mathfrak{gl}_N)$.
\end{proposition}

\begin{definition}[\cite{MRS}]
\label{Uqon}
The twisted quantum loop algebra $\mathfrak{U}_q(\mathcal{L}^{tw}(\mathfrak{o}_N))$ is the unital associative algebra
over $\mathbb{C}(q)$ generated by $\{S_{ij}^{(r)} \,|\, r \in \mathbb{Z}_+, \, 1 \leq i, j \leq N\}$ which are subject to the
relations
\begin{gather*}
S_{ij}^{(0)}=0 \qquad \text{if} \quad  1 \leq i < j \leq n,
\qquad
S_{ii}^{(0)}=1 \qquad \forall\,  1 \leq i \leq n,
\\
R_q(u,v)S_1(u)R_q^t\big(u^{-1},v\big)S_2(v)=S_2(v)R_q^t\big(u^{-1},v\big)S_1(u)R_q(u,v),
\end{gather*}
where $S_a(u)$ is obtained from $S(u)$ in the same way as in Def\/inition~\ref{Ytw}, except that in this case we have
$S(u):=\sum\limits_{i,j=1}^N E_{ij} \otimes S_{ij}(u)$ and $S_{ij}(u):=\sum\limits_{r=0}^\infty S_{ij}^{(r)}u^{-r}$.
\end{definition}

\begin{definition}[\cite{MRS}]
%\label{Uqspn}
The twisted quantum loop algebra $\mathfrak{U}_q(\mathcal{L}^{tw}(\mathfrak{sp}_N))$ is the unital associative algebra
over $\mathbb{C}(q)$ with generators $\big\{S_{ij}^{(r)} \,|\, r \in \mathbb{Z}_+, 1 \leq i, j \leq N\big\}$ and
$\big\{(S_{ii'}^{(0)})^{-1} \,|\, i=1,3,\ldots,N{-}1\big\}$, which are subject to the relations
\begin{gather*}
S_{ij}^{(0)}=0 \qquad \text{whenever} \quad i<j \quad \text{and} \quad j \neq i'
\\
S_{i'i'}^{(0)}S_{ii}^{(0)}-q^2S_{i'i}^{(0)}S_{ii'}^{(0)}=q^3, \qquad i=1, 3, \ldots, N-1,
\\
S_{ii'}^{(0)}(S_{ii'}^{(0)})^{-1}=(S_{ii'}^{(0)})^{-1}S_{ii'}^{(0)}=1, \qquad i=1, 3, \ldots, N-1,
\\
R_q(u,v)S_1(u)R_q^t\big(u^{-1},v\big)S_2(v)=S_2(v)R_q^t\big(u^{-1},v\big)S_1(u)R_q(u,v),
\end{gather*}
where $S_a(u)$ is def\/ined here in the same way as in Def\/inition~\ref{Uqon}.
\end{definition}

The twisted quantum loop algebra $\mathfrak{U}_q(\mathcal{L}^{tw}(\mathfrak{g}_N))$ is a~deformation of the universal
enveloping algebra of the twisted loop algebra $\mathfrak{g}_N^{tw}[s,s^{-1}]$ which is def\/ined in the following way:

\begin{definition}
The twisted loop algebra $\mathfrak{g}_N^{tw}[s,s^{-1}]$ is the Lie subalgebra of
$\mathcal{L}(\mathfrak{gl}_N)$ given~by
\begin{gather*}
\mathfrak{g}_N^{tw}\big[s,s^{-1}\big]=\big\{A(s) \in \mathcal{L}(\mathfrak{gl}_N) \,|\, \sigma(A(s))=A\big(s^{-1}\big)\big\},
\end{gather*}
where~$\sigma$ is given by~\eqref{orthauto} for $\mathfrak{g}_N=\mathfrak{o}_N$ and by~\eqref{sympauto} for
$\mathfrak{g}_N=\mathfrak{sp}_N$.
This algebra is also denoted by $\mathcal{L}^{tw}(\mathfrak{g}_N)$.
\end{definition}

The twisted quantum loop algebras can be regarded as subalgebras of
$\mathfrak{U}_q(\mathcal{L}(\mathfrak{gl}_N))$:

\begin{proposition}[\cite{MRS}]\label{embU}
The assignment
\begin{gather*}
S_{ij}^{(r)} \mapsto \sum\limits_{k,l=1}^N \sum\limits_{p=0}^r b_{kl} T_{ik}^{(p)} \overline{T}_{jl}^{(r-p)}
\end{gather*}
provides an embedding of $\mathfrak{U}_q(\mathcal{L}^{tw}(\mathfrak{g}_N))$ into
$\mathfrak{U}_q(\mathcal{L}(\mathfrak{gl}_N))$.
\end{proposition}

For each $r>0$, let $S_{ij}^{(r,0)}=\frac{S_{ij}^{(r)}}{q-q^{-1}}$.
When $\mathfrak{g}_N = \mathfrak{o}_N$, we set $S_{ij}^{(0,0)}=\frac{S_{ij}^{(0)}}{q-q^{-1}}$ when $i>j$ and
$S_{ij}^{(0,0)}=-S_{ji}^{(0,0)}$ when $i\le j$; when $\mathfrak{g}_N = \mathfrak{sp}_N$, we set
$S_{ij}^{(0,0)}=\frac{S_{ij}^{(0)}-b_{ij}}{q-q^{-1}}$ when $i\ge j$ or $j=i'$ and $S_{ij}^{(0,0)} = -S_{ji}^{(0,0)}$
when $i<j$ and $j\neq i'$.

We def\/ine inductively $S_{ij}^{(r,m)}$ by
\begin{gather*}
S_{ij}^{(r,m+1)}=S_{ij}^{(r+1,m)}-S_{ij}^{(r,m)}
\end{gather*}
for $m \geq 0$.
\begin{lemma}
Identify $S_{ij}^{(r)}$ with its image under the embeddings of Proposition~{\rm \ref{embU}}.
For any $r>0$ and $m \geq 0$, we have
\begin{gather*}
S_{ij}^{(r,m)}=\sum\limits_{k=1}^N \left(b_{kj} T_{ik}^{(r,m)} + b_{ik} \overline{T}_{jk}^{(r,m)} \right) +
\big(q-q^{-1}\big) \sum\limits_{k,l=1}^N \sum\limits_{p=1}^{r-1} b_{kl} T_{ik}^{(p,0)} \overline{T}_{jl}^{(r-p,m)}
\\
\phantom{S_{ij}^{(r,m)}=}{}
 + \big(q-q^{-1}\big)\!\left(\sum\limits_{\substack{1 \leq k \leq i\\ 1 \leq l \leq N}}\!\left(\frac{q}{q+1}\right)^{\delta_{ik}}\!
b_{kl} T_{ik}^{(0,0)} \overline{T}_{jl}^{(r,m)} + \sum\limits_{\substack{1 \leq k \leq N \\ j \leq l \leq N}}\!
\left(\frac{q}{q+1}\right)^{\delta_{jl}} \! b_{kl} T_{ik}^{(r,m)} \overline{T}_{jl}^{(0,0)} \!\right)
\\
\phantom{S_{ij}^{(r,m)}=}{}
  + \big(q-q^{-1}\big) \sum\limits_{k,l=1}^N \sum\limits_{a+b=m-1} b_{kl} T_{ik}^{(r,a)} \overline{T}_{jl}^{(1,b)}.
\end{gather*}
\end{lemma}

\begin{proof}
Since $T_{ik}^{(0)}=0$ when $i<k$ and $\overline{T}_{jl}^{(0)}=0$ when $j>l$, we have by def\/inition
\begin{gather*}
S_{ij}^{(r,0)}= \big(q-q^{-1}\big) \sum\limits_{k,l=1}^N \sum\limits_{p=1}^{r-1} b_{kl} \left(\frac{T_{ik}^{(p)}}{q-q^{-1}}
\right) \left(\frac{\overline{T}_{jl}^{(r-p)}}{q-q^{-1}} \right)
\\
\phantom{S_{ij}^{(r,0)}=}{}
+ \big(q-q^{-1}\big)\sum\limits_{\substack{1 \leq k < i\\ 1 \leq l \leq N}} b_{kl} \left(\frac{T_{ik}^{(0)}}{q-q^{-1}} \right)
\left(\frac{\overline{T}_{jl}^{(r)}}{q-q^{-1}}\right)
\\
\phantom{S_{ij}^{(r,0)}=}{}
+ (q-1) \sum\limits_{l=1}^N b_{il} \left(\frac{T_{ii}^{(0)}-1}{q-1} \right)
\left(\frac{\overline{T}_{jl}^{(r)}}{q-q^{-1}}\right) + \sum\limits_{l=1}^N b_{il}
\frac{\overline{T}_{jl}^{(r)}}{q-q^{-1}}
\\
\phantom{S_{ij}^{(r,0)}=}{}
+ \big(q-q^{-1}\big) \sum\limits_{\substack{1 \leq k \leq N\\ j < l \leq N}} b_{kl} \left(\frac{T_{ik}^{(r)}}{q-q^{-1}} \right)
\left(\frac{\overline{T}_{jl}^{(0)}}{q-q^{-1}}\right)
\\
\phantom{S_{ij}^{(r,0)}=}{}
+ (q-1) \sum\limits_{k=1}^N b_{kj} \left(\frac{T_{ik}^{(r)}}{q-q^{-1}} \right)
\left(\frac{\overline{T}_{jj}^{(0)}-1}{q-1}\right) + \sum\limits_{k=1}^N b_{kj} \frac{T_{ik}^{(r)}}{q-q^{-1}}
\\
\phantom{S_{ij}^{(r,0)}}{}
= \sum\limits_{k=1}^N \left(b_{kj} T_{ik}^{(r,0)} + b_{ik} \overline{T}_{jk}^{(r,0)} \right) + \big(q-q^{-1}\big)
\sum\limits_{k,l=1}^N \sum\limits_{p=1}^{r-1} b_{kl} T_{ik}^{(p,0)} \overline{T}_{jl}^{(r-p,0)}
\\
\phantom{S_{ij}^{(r,0)}=}{}
+ \big(q-q^{-1}\big) \! \left(\sum\limits_{\substack{1 \leq k \leq i\\ 1 \leq l \leq N}}\! \left(\frac{q}{q+1}\right)^{\delta_{ik}}\!
b_{kl} T_{ik}^{(0,0)} \overline{T}_{jl}^{(r,0)} + \sum\limits_{\substack{1 \leq i \leq N\\ j \leq l \leq N}}\!
\left(\frac{q}{q+1}\right)^{\delta_{jl}} \! b_{kl} T_{ik}^{(r,0)} \overline{T}_{jl}^{(0,0)} \right).
\end{gather*}
This proves the case $m=0$.
The case $m=1$ is similar, except for the presence of an extra sum.
The general case follows immediately by induction on~$m$.
\end{proof}

The previous lemma along with~\eqref{TbarT} yields the next corollary.

\begin{corollary}
\label{Scong}
Under the same assumption as in the previous lemma, we have $S_{ij}^{(r,m)} \in \mathbf{K}_m$, and
\begin{gather*}
S_{ij}^{(r,m)} \equiv \sum\limits_{k=1}^N \left(b_{kj}T_{ik}^{(0,m)} + (-1)^{m+1} b_{ik} T_{jk}^{(0,m)}\right) \\
\hphantom{S_{ij}^{(r,m)} \equiv}{}
+
\big(q-q^{-1}\big) \sum\limits_{k,l=1}^N \sum\limits_{p=1}^m (-1)^{m+1-p} b_{kl} T_{ik}^{(0,p-1)} T_{jl}^{(0,m-p)}
\end{gather*}
modulo $\mathbf{K}_{m+1}$.
\end{corollary}
Again let's view $S_{ij}^{(r)}$ as an element of $\mathfrak{U}_q(\mathcal{L}(\mathfrak{gl}_N))$.
Let $\zeta_{ij}^{(r,m)}$ be the image of $S_{ij}^{(r,m)}$ in $\mathbf{K}_m/\mathbf{K}_{m+1}$.
Note that, by Corollary~\ref{Scong}, $\zeta_{ij}^{(r,m)}$ is independent of $r>0$.
\begin{theorem}
\label{twUqY}
Identify $s_{ij}^{(m+1)}$ with its image under the embedding of Proposition~{\rm \ref{embY}}.
For each $m \geq 0$, we have
\begin{gather*}
\varphi\big(s_{ij}^{(m+1)}\big)=\zeta_{ij}^{(1,m)},
\end{gather*}
where~$\varphi$ is the isomorphism of Theorem~{\rm \ref{mainthm}}.
\end{theorem}

\begin{remark}
It would also be possible to prove a~similar theorem with $\zeta_{ij}^{(0,m)}$ instead of $\zeta_{ij}^{(1,m)}$, but the
proof would be longer because we would have to consider dif\/ferent cases in order to take into account that about half of
the generators $S_{ij}^{(0)}$ are~$0$.
\end{remark}

\begin{proof}
By Proposition~\ref{embY}, Corollary~\ref{Scong} and the def\/inition of~$\varphi$, we have
\begin{gather*}
\varphi\big(s_{ij}^{(m+1)}\big)=\sum\limits_{k=1}^N \left(g_{kj} \xi_{ik}^{(0,m)} + (-1)^{m+1} g_{ik} \xi_{jk}^{(0,m)}\right)
\\
\hphantom{\varphi\big(s_{ij}^{(m+1)}\big)=}{}
+ \hbar \sum\limits_{k,l=1}^N \sum\limits_{p=1}^m (-1)^{m+1-p} g_{kl} \xi_{ik}^{(0,p-1)} \xi_{jl}^{(0,m-p)}
  = \zeta_{ij}^{(1,m)},
\end{gather*}
where we have used the fact that $b_{kl}-g_{kl}$ is divisible by $q-1$.
\end{proof}

We can obtain an analogue of Theorem~\ref{mainthm}.
From Theorem~\ref{qlim} and Proposition~\ref{embU}, we can deduce that the enveloping algebra of
$\mathfrak{g}_N^{tw}[s]$ is the limit when $q\rightarrow 1$ of $\mathfrak{U}_q(\mathcal{L}^{tw}(\mathfrak{g}_N))$ in the
following sense: if we let $\mathfrak{U}_{\mathcal{A}}(\mathcal{L}^{tw}(\mathfrak{g}_N))$ be the
$\mathcal{A}$-subalgebra of $\mathfrak{U}_q(\mathcal{L}^{tw}(\mathfrak{g}_N))$ generated by the $S_{ij}^{(r,0)}$ for all
$i$, $j$ and all $r\ge 0$, then
$\mathfrak{U}_{\mathcal{A}}(\mathcal{L}^{tw}(\mathfrak{g}_N))/(q-1)\mathfrak{U}_{\mathcal{A}}(\mathcal{L}^{tw}(\mathfrak{g}_N))$
is isomorphic to $\mathfrak{U}(\mathfrak{g}_N^{tw}[s,s^{-1}])$ (see the proof of Corollaries~3.5 and~3.12
in~\cite{MRS}).

We can def\/ine an algebra $\widetilde{Y}^{tw}(\mathfrak{g}_N)$ similarly to how we def\/ined
$\widetilde{Y}(\mathfrak{gl}_N)$.
For $m\ge 0$, in the orthogonal case, denote by $\mathsf{K}_m^{tw}$ the Lie ideal of $\mathfrak{o}_N^{tw}[s,s^{-1}]$
spanned~by
\begin{gather*}
E_{ij} s^r(s-1)^m - E_{ji} s^{-r} \big(s^{-1}-1\big)^m
\end{gather*}
for all $r \in \mathbb{Z}$.
In the symplectic case, let $\mathsf{K}_m^{tw}$ be the Lie ideal of $\mathfrak{sp}_N^{tw}[s,s^{-1}]$ spanned~by
\begin{gather*}
E_{ij'}s^r(s-1)^m - (-1)^{i+j+1} E_{ji'} s^{-r}\big(s^{-1}-1\big)^m
\end{gather*}
for all $r \in \mathbb{Z}$.
Let $U^{tw}$ be the subspace of $\mathfrak{U}_{\mathcal{A}}(\mathcal{L}^{tw}(\mathfrak{g}_N))$ spanned over $\mathbb{C}$
by all the genera\-tors~$S_{ij}^{(r,0)}$, and observe that $U^{tw}\cap
(q-1)\mathfrak{U}_{\mathcal{A}}(\mathcal{L}^{tw}(\mathfrak{g}_N)) = \{0 \}$.
Let $\mathbb{K}_m^{tw}$ be the two-sided ideal of $\mathfrak{U}_{\mathcal{A}}(\mathcal{L}^{tw}(\mathfrak{g}_N))$
generated by $\psi^{-1}(\mathsf{K}_m^{tw})\cap U^{tw}$, where~$\psi$ this time denotes the composite
\begin{gather*}
\mathfrak{U}_{\mathcal{A}}(\mathcal{L}^{tw}(\mathfrak{g}_N)) \twoheadrightarrow
\mathfrak{U}_{\mathcal{A}}(\mathcal{L}^{tw}(\mathfrak{g}_N))/(q-1)\mathfrak{U}_{\mathcal{A}}(\mathcal{L}^{tw}(\mathfrak{g}_N))
\stackrel{\sim}{\longrightarrow} \mathfrak{U}\big(\mathfrak{g}_N^{tw}\big[s,s^{-1}\big]\big).
\end{gather*}
Set $\mathbf{K}_m^{tw}$ equal to the sum of the ideals $(q-q^{-1})^{m_0} \mathbb{K}_{m_1}^{tw}\cdots
\mathbb{K}_{m_k}^{tw}$ with $m_0+m_1+\cdots + m_k \geq m$.

Let $\widetilde{Y}^{tw}(\mathfrak{g}_N)$ be the $\mathbb{C}$-algebra
\begin{gather*}
\bigoplus_{m=0}^{\infty} \mathbf{K}_m^{tw}/\mathbf{K}_{m+1}^{tw},
\end{gather*}
where $\mathbf{K}_0^{tw} = \mathfrak{U}_{\mathcal{A}}(\mathcal{L}^{tw}(\mathfrak{g}_N))$.
We also view $\widetilde{Y}^{tw}(\mathfrak{g}_N)$ as an algebra over $\mathbb{C}[\hbar]$ by setting $\hbar =
\overline{q-q^{-1}} \in \mathbf{K}_1^{tw}/\mathbf{K}_2^{tw}$.

\begin{theorem}
\label{twUqY2}
$Y^{tw}(\mathfrak{g}_N)$ is isomorphic to $\widetilde{Y}^{tw}(\mathfrak{g}_N)$ via the function $\varphi^{tw}$ that
sends $s_{ij}^{(m+1)}$ to $\overline{S_{ij}^{(1,m)}} \in \mathbf{K}_{m}^{tw}/\mathbf{K}_{m+1}^{tw}$ for $m \ge 0$.
\end{theorem}

\begin{proof}
Theorem~\ref{twUqY} implies that the following diagram is commutative:
\begin{gather*}
\xymatrix{Y^{tw}(\mathfrak{g}_N) \ar[r] \ar[d]_{\varphi^{tw}} & Y(\mathfrak{gl}_N) \ar[d]^{\varphi}\\
 \widetilde{Y}^{tw}(\mathfrak{g}_N) \ar[r] & \widetilde{Y}(\mathfrak{gl}_N).}
\end{gather*}
In this diagram, the top horizontal arrow is the embedding of Proposition~\ref{embY} and the bottom horizontal arrow is
the one induced from the embedding of Proposition~\ref{embU}.
The injectivity of~$\varphi^{tw}$ now follows from the fact that~$\varphi$ provides an isomorphism between
$Y(\mathfrak{gl}_N)$ and $\widetilde{Y}(\mathfrak{gl}_N)$: see Theorem~\ref{mainthm}.

We need to see that $\varphi^{tw}$ is surjective.
Def\/ine elements $\widetilde{S}_{ij}^{(r,m)}$ with $0 \leq r \leq m$ as follows: for each $m \geq 0$, let
$\widetilde{S}_{ij}^{(0,m)}=S_{ji}^{(0,m)}$ and $\widetilde{S}_{ij}^{(m,m)}=(-1)^{m+1}S_{ij}^{(0,m)}$, and for $1 \leq r
\leq m$ let
\begin{gather*}
\widetilde{S}_{ij}^{(r,m+1)}=\widetilde{S}_{ij}^{(r-1,m)}-\widetilde{S}_{ij}^{(r,m)}.
\end{gather*}

Then induction on~$m$ shows that
\begin{gather*}
\psi(S_{ij}^{(r,m)})=E_{ij}s^r(s-1)^m - E_{ji}s^{-r}\big(s^{-1}-1\big)^m,\\
\psi(S_{ji}^{(r,m)})=(-1)^{m+1}\big(E_{ij}s^{-(m+r)}(s-1)^m - E_{ji}s^{m+r}\big(s^{-1}-1\big)^m\big),\\
\psi(\widetilde{S}_{ij}^{(r,m)})=(-1)^{m+1}\big(E_{ij}s^{-(m-r)}(s-1)^m - E_{ji}s^{m-r}\big(s^{-1}-1\big)^m\big).
\end{gather*}
in the orthogonal case, and similarly in the symplectic case.
It follows that for any f\/ixed~$m$, the images of these elements under~$\psi$ span $\mathsf{K}_m^{tw}$, and they are all
in $U^{tw}$ by def\/inition.
Now note that for any element $X \in \psi^{-1}(\mathsf{K}_m^{tw}) \cap U^{tw}$, there is some element~$Y$ in
 \begin{gather*}
 \operatorname{span}_{\mathbb{C}} \big\{S_{ij}^{(r,m)},\widetilde{S}_{ij}^{(r,m)} \,|\, i,j=1,\ldots,N,\; r \in \mathbb{Z}_{\ge 0}\big\}
 \end{gather*}
for which $X-Y \in (q-1)\mathfrak{U}_{\mathcal{A}}(\mathcal{L}^{tw}(\mathfrak{g}_N))$.
This follows from the fact that
\begin{gather*}
\mathfrak{U}_{\mathcal{A}}\big(\mathcal{L}^{tw}(\mathfrak{g}_N)\big)/(q-1)\mathfrak{U}_{\mathcal{A}}\big(\mathcal{L}^{tw}(\mathfrak{g}_N)\big)
\stackrel{\sim}{\longrightarrow} \mathfrak{U}\big(\mathfrak{g}_N^{tw}\big[s,s^{-1}\big]\big)
\end{gather*}
is an isomorphism.
Since $X-Y$ is also in $U^{tw}$ and since $U^{tw} \cap
(q-1)\mathfrak{U}_{\mathcal{A}}(\mathcal{L}^{tw}(\mathfrak{g}_N))=\{0\}$, we see that $X=Y$.
Therefore, any element of $\mathbf{K}_m^{tw}$ is a~sum of monomials $f(q)(q-q^{-1})^{m_0}\mathcal{M}$, where $f(q) \in
\mathcal{A}$ is not divisible by $q-1$ and $\mathcal{M}=\sigma_{i_1 j_1}^{(r_1,m_1)} \cdots \sigma_{i_k j_k}^{(r_k,m_k)}$ with
\begin{gather*}
\sigma_{i_d j_d}^{(r_d,m_d)} \in \big\{S_{i_d j_d}^{(r_d,m_d)},\widetilde{S}_{i_d j_d}^{(r_d,m_d)}\big\}
\end{gather*}
and $m_0+\cdots+m_k \geq m$.
Following the same argument as in the $\mathfrak{gl}_N$ case, such a~monomial is congruent modulo
$\mathbf{K}_{m+1}^{tw}$ to
\begin{gather*}
f(1)\big(q-q^{-1}\big)^{m_0}S_{i_1 j_1}^{(1,m_1)} \cdots S_{i_k j_k}^{(1,m_k)}
\end{gather*}
up to a~sign.
The image modulo $\mathbf{K}_{m+1}^{tw}$ of this element is
\begin{gather*}
f(1)\hbar^{m_0}\overline{S_{i_1 j_1}^{(1,m_1)}} \cdots \overline{S_{i_k j_k}^{(1,m_k)}}
\end{gather*}
and this is in the image of $\varphi^{tw}$, which proves that $\varphi^{tw}$ is surjective.
\end{proof}

\subsection*{Acknowledgements}

The research of N.~Guay was supported by a~Discovery Grant from the Natural Sciences and Engineering Research Council of
Canada.
P.~Conner was supported by a~Postgraduate Scholarship CGS-M from the same agency.

\pdfbookmark[1]{References}{ref}
\LastPageEnding

\end{document}